\documentclass{amsart}
\usepackage{amssymb}

\title{A Variant of $K$-theory: $K_{\pm }$}

\author{Michael~Atiyah}
\address{School of Mathematics 
\\  University of Edinburgh \\ James Clerk Maxwell Buildings
\\ King's Buildings
\\ Edinburgh \ EH9 3JZ}

\author{Michael~Hopkins}
\address{Department of Mathematics \\ Massachusetts Institute of
Technology\\Cambridge, MA 02139-4307}
\email{mjh@math.mit.edu}

\DeclareMathOperator{\coker}{coker}

\theoremstyle{remark}
	\newtheorem{rem}{Remark}

\begin{document}

\maketitle

\section*{Introduction}

Topological $K$-theory~\cite{atiyah67:_k} has many variants which have
been developed and exploited for geometric purposes. There are real or
quaternionic versions, ``Real'' $K$-theory in the sense
of~\cite{atiyah66:_k}, equivariant $K$-theory~\cite{segal68:_equiv_k}
and combinations of all these.

In recent years $K$-theory has found unexpected application in the
physics of string theories~\cite{e_gauge_theor_k_theor_m_theor,
minasian97:_k_ramon_ramon,moore00:_self_ramon_ramon_k, witten98:_d_k}
and all variants of $K$-theory that had previously been developed
appear to be needed. There are even variants, needed for the physics,
which had previously escaped attention, and it is one such variant
that is the subject of this paper.

This variant, denoted by $K_{\pm }(X),$ was introduced by
Witten~\cite{witten98:_d_k} in relation to ``orientifolds''. The
geometric situation concerns a manifold $X$ with an involution
$\tau $ leaving a fixed sub-manifold $Y$.   On $X$ one wants to
study a pair of complex vector bundles $(E^{+},E^{-})$ with the
property that $\tau$ interchanges them. If we think of the virtual
vector bundle $E^{+}-E^{-}$, then $\tau $ takes this into its
negative, and $K_{\pm }(X)$ is meant to be the appropriate $K$-theory
of this situation.

In physics, $X$ is a 10-dimensional Lorentzian manifold and maps $\Sigma
\rightarrow X$ of a surface $\Sigma $ describe the world-sheet of strings.
The symmetry requirements for the appropriate Feynman integral impose
conditions that the putative $K$-theory $K_{\pm }(X)$ has to satisfy.

The second author proposed a precise topological definition of $K_{\pm
}(X)$ which appears to meet the physics requirements, but it was not
entirely clearly how to link the physics with the geometry.

In this paper we elaborate on this definition and also a second (but
equivalent) definition of $K_{\pm }(X)$. Hopefully this will bring the
geometry and physics closer together, and in particular link it up with the
analysis of Dirac operators.

Although $K_{\pm}(X)$ is defined in the context of spaces with
involution it is rather different from Real $K$-theory or equivariant
$K$-theory (for $G=Z_{2}$) although it has superficial resemblances
to both. The differences will become clear as we proceed but at this
stage it may be helpful to consider the analogy with
cohomology. Equivariant cohomology can be defined (for any compact Lie
group $G$), and this has relations with equivariant $K$-theory. But
there is also ``cohomology with local coefficients'', where the
fundamental group $\pi_{1}(X)$ acts on the abelian coefficient
group. In particular for integer coefficients $Z$ the only such
action is via a homomorphism $\pi_{1}(X)\rightarrow Z_{2}$, i.e. by
an element of $H^{1}(X;Z_{2})$ or equivalently a double-covering
$\tilde{X}$ of $X$.

This is familiar for an unoriented manifold with $\tilde{X}$ its oriented
double-cover. In this situation, if say $X$ is a compact $n$-dimensional
manifold, then we do not have a fundamental class in $H^{n}(X;Z)$ but in $
H^{n}(X;\tilde{Z})$ where $\tilde{Z}$ is the local coefficient system
defined by $\tilde{X}$. This is also called ``twisted cohomology''.

Here $\tilde{X}$ has a fixed-point-free involution $\tau $ and, in such a
situation, our group $K_{\pm }(\tilde{X})$ is the precise $K$-theory
analogue of twisted cohomology. This will become clear later.

In fact $K$-theory has more sophisticated twisted versions.
In~\cite{donovan70:_graded_brauer_k} Donovan and Karoubi use Wall's
graded Brauer group~\cite{wall63:_graded_brauer} to construct
twistings from elements of $H^{1}(X;Z_{2})\times
H^{3}(X;Z)_{\text{torsion}}$.  More general twistings of $K$-theory
arise from automorphisms of its classifying space, as do twistings of
equivariant $K$-theory.  Among these are twistings involving a general
element of $H^{3}(X;Z)$ (i.e., one which is not necessarily of finite
order).  These are also of interest in physics, and have recently been
the subject of much attention~\cite{atiyah:_twist_k,
bouwknegt02:_twist_k_k,math.AT/0206257}.  Our $K_{\pm}$ is a twisted
version of equivariant $K$-theory,\footnote{It is the twisting of
equivariant $K$-theory by the non-trivial element of
$H^{1}_{Z_{2}}(\text{pt})=Z_{2}$.  From the point of view of the
equivariant graded Brauer group, $K_{\pm}(X)$ is the $K$-theory of the
graded cross product algebra $C(X)\otimes M\rtimes Z_{2}$, where $C(X)$
is the algebra of continuous functions on $X$ and $M$ is the graded
algebra of $2\times 2$-matrices over the complex numbers, graded in
such a way that $(i,j)$ entry has degree $i+j$. The action of $Z_{2}$ is
the combination of the geometric action given on $X$ and conjugation
by the permutation matrix on $M$.} and this paper can be seen as a
preliminary step towards these other more elaborate version.

\section*{2.  The first definition}

Given a space $X$ with involution we have two natural $K$-theories,
namely $K(X)$ and $K_{Z_{2}}(X)$ --- the ordinary and equivariant
theories respectively. Moreover we have the obvious homomorphism
\begin{equation}
\phi :K_{Z_{2}}(X)\rightarrow K(X)  \tag{2.1}
\end{equation}
which ``forgets'' about the $Z_{2}$-action. We can reformulate this by
introducing the space $(X\times Z_{2})$ with the involution $
(x,0)\rightarrow (\tau (x),1)$. Since this action is free we have 
\[
K_{Z_{2}}(X\times Z_{2})\cong K(X) 
\]
and (2.1)\ can then be viewed as the natural homomorphism for $K_{Z_{2}}$ 
induced by the projection 
\begin{equation}
\pi :X\times Z_{2}\rightarrow X\;.  \tag{2.2}
\end{equation}

Now, whenever we have such a homomorphism, it is part of a long exact
sequence (of period 2) which we can write as an exact triangle 
\begin{equation}
\begin{tabular}{ccc}
$K_{Z_{2}}^{\ast }(X)$ & $\stackrel{\phi }{\rightarrow }$ & $K^{\ast }(X)$
\\ 
$\nwarrow $ &  & $\swarrow \delta $ \\ 
& $K_{Z_{2}}^{\ast }(\pi )$ & 
\end{tabular}
\tag{2.3}
\end{equation}
where $K^{\ast }=K^{0}\oplus K^{1},\;\delta$ has degree $1\mod 2$
and the relative group $K_{Z_{2}}^{\ast }(\pi )$ is just the relative
group for a pair, when we replace $\pi $\ by a $Z_{2}$-homotopically
equivalent inclusion. In this case a natural way to do this is to
replace the $X$ factor on the right of (2.2) by $X\times I$ where
$I=[0,1]$ is the unit interval with $\tau$ being reflection about the
mid-point $\frac{1}{2}$. Thus, explicitly,
\begin{equation}
K_{Z_{2}}^{\ast }(\pi )=K_{Z_{2}}^{\ast }(X\times I,X\times \partial I) 
\tag{2.4}
\end{equation}
where $\partial I$ is the (2-point) boundary of $I.$

We now take the group in (2.4) (with the degree shifted by one) as our
definition of $K_{\pm }^{\ast }(X).$ It is then convenient to follow the
notation of~\cite{atiyah66:_k} where $R^{p,q}=R^{p}\oplus R^{q}$ with the involution
changing the sign of the first factor, and we use $K$-theory with compact
supports (so as to avoid always writing the boundary). Then our definition
of $K_{\pm }$ becomes 
\begin{equation}
K_{\pm }^{0}(X)=K_{Z_{2}}^{1}(X\times R^{1,0})\cong K_{Z_{2}}^{0}(X\times
R^{1,1})  \tag{2.5}
\end{equation}
(and similarly for $K^{1}).$

Let us now explain why this definition fits the geometric situation we began
with (and which comes from the physics). Given a vector bundle $E$ we can
form the pair $(E,\tau ^{\ast }E)$ or the virtual bundle 
\[
E-\tau ^{\ast }E. 
\]
Under the involution, $E$ and $\tau ^{\ast }E$ switch and the virtual
bundle goes into its negative. Clearly, if $E$ came from an equivariant
bundle, then $E\cong \tau ^{\ast }E$ and the virtual bundle is zero. Hence
the virtual bundle depends only the element defined by $E$ in the cokernel
of $\phi $, and hence by the image of $E$ in the next term of the exact
sequence (2.3) i.e. by 
\[
\delta (E)\in K_{\pm }^{0}(X). 
\]

This explains the link with our starting point and it also shows that one
cannot always define $K_{\pm }(X)$ in terms of such virtual bundles on $
X.\; $In general the exact sequence $(2.3)$ does not break up into short
exact sequences and $\delta $ is not surjective.

At this point a physicist might wonder whether the definition of $K_{\pm
}(X)$ that we have given is the right one. Perhaps there is another group
which is represented by virtual bundles. We will give two pieces of evidence
in favour of our definition, the first pragmatic and the second more
philosophical.

First let us consider the case when the involution $\tau $ on $X$ is
trivial. Then

\noindent $K_{Z_{2}}^{\ast }\left( X\right) =R(Z_{2})\otimes K^{\ast }(X)$ 
and $R(Z_{2})=Z\oplus Z$ is the representation ring of $Z_{2}$ and is
generated by the two representations 
\[
\begin{tabular}{ll}
1 & (trivial representation) \\ 
$\rho $ & (sign representation).
\end{tabular}
\]
The homomorphism $\phi $ is surjective with kernel $(1-\rho )K^{\ast }(X)$ 
so $\delta =0$ and  
\begin{equation}
K_{\pm }^{0}(X)\cong K^{1}(X)\;.  \tag{2.6}
\end{equation}
This fits with the requirements of the physics, which involves a switch from
type $\mathrm{II}A$ to type $\mathrm{II}B$\ string theory. Note also that
it gives an extreme example when $\partial $ is not surjective.

Our second argument is concerned with the general passage from physical
(quantum) theories to topology. If we have a theory with some symmetry then
we can consider the quotient theory, on factoring out the symmetry.
Invariant states of the original theory become states of the quotient theory
but there may also be new states that have to be added. For example if we
have a group $G$ of geometric symmetries, then closed strings in the
quotient theory include strings that begin at a point $x$ and end at $
g(x)\; $for $g\in G.$ All this is similar to what happens in topology with
(generalized) cohomology theories, such as $K$-theory. If we have a morphism
of theories, such as $\phi $ in\ (2.1)\ then the third theory we get fits
into a long-exact sequence. The part coming from $K(X)$ is only part of the
answer -- other elements have to be added. In ordinary cohomology where we
start with cochain complexes the process of forming a quotient theory
involves an ordinary quotient (or short exact sequence) at the level of
cochains. But this becomes a long exact sequence at the cohomology level.
For $K$-theory the analogue is to start with bundles over small open sets
and at this level we can form the naive quotients, but the $K$-groups arise
when we impose the matching conditions to get bundles, and then we end up
with long exact sequences.

It\ is also instructive to consider the special case when the involution is
free so that we have a double covering $\tilde{X}\rightarrow X$ and the
exact triangle (2.3), with $\tilde{X}$ for $X$, becomes the exact triangle 
\begin{equation}
\begin{tabular}{ccc}
$K^{\ast }(X)$ & $\stackrel{\phi }{\rightarrow }$ & $K^{\ast }(\tilde{X})$
\\ 
$\nwarrow $ &  & $\swarrow \delta $ \\ 
& $K_{Z_{2}}^{\ast }(L)$ & 
\end{tabular}
\tag{2.7}
\end{equation}
Here $L$ is the real line bundle over $X$ associated to the double
covering $\tilde{X}$ (or to the corresponding element of $
H^{1}(X,Z_{2})),\; $and we again use compact supports. Thus (for
$q=0,1 \mod 2)$
\begin{equation}
K_{\pm }^{q}(\tilde{X})=K^{q+1}(L).  \tag{2.8}
\end{equation}

If we had repeated this argument using equivariant cohomology instead of
equivariant $K$-theory we would have ended up with the twisted cohomology
mentioned earlier, via a twisted suspension isomorphism, 
\begin{equation}
H^{q}(X,\tilde{Z})=H^{q+1}(L).  \tag{2.9}
\end{equation}
This shows that, for free involutions, $K_{\pm }$ \textbf{is precisely the
analogue of twisted cohomology,\ }so that, for example, the Chern character
of the former take values in the rational extension of the latter.

\section*{3. Relation to Fredholm operators}

In this section we shall given another definition of $K_{\pm }$ which ties
in naturally with the analysis of Fredholm operators, and we shall show that
this definition is equivalent to the one given in Section 2.

We begin by recalling a few basic facts about Fredholm operators. Let $H$ 
be complex Hilbert space, $\mathcal{B\;}$the space of bounded operators with
the norm topology and $\mathcal{F\subset B\;}$the open subspace of Fredholm
operators, i.e. operators $A$ so that $\ker A$ and $\coker A$ are
both finite-dimensional. The index defined by 
\[
\text{index }A=\dim \ker A-\dim \coker A 
\]
is then constant on connected components of $\mathcal{F}.$ If we introduce
the adjoint $A^{\ast }$ of $A$ then 
\[
\coker A=\ker A^{\ast } 
\]
so that 
\[
\text{index }A=\dim \ker A-\dim \ker A^{\ast }. 
\]

More generally if we have a continuous map 
\[
f:X\rightarrow \mathcal{F} 
\]
(i.e. a family of Fredholm operators, parametrized by $X$), then one can
define 
\[
\text{index }f\in K(X) 
\]
and one can show~\cite{atiyah67:_k} that we have an isomorphism 
\begin{equation}
\text{index}:[X,\mathcal{F]}\cong K(X)  \tag{3.1}
\end{equation}
where $[\phantom{x},\phantom{x}]$ denotes homotopy classes of
maps. Thus $K(X)$ has a natural definition as the ``home'' of indices
of Fredholm operators (parametrized by $X)$:\ it gives the complete
homotopy invariant.

Different variants of $K$-theory can be defined by different variants of
(3.1). For example real $K$-theory uses real Hilbert space and equivariant $
K $-theory for $G$-spaces uses a suitable $H$-space module of $G$, namely $
L_{2}(G)\otimes H.$ It is natural to look for a similar story for our new
groups $K_{\pm }(X).$ A first candidate might be to consider $Z_{2}$
-equivariant maps 
\[
f:X\rightarrow \mathcal{F} 
\]
where we endow $\mathcal{F\;}$with the involution $A\rightarrow A^{\ast }$ 
given by taking the adjoint operator. Since this switches the role of kernel
and cokernel it acts as $-1$ on the index, and so is in keeping with our
starting point.

As a check we can consider $X$ with a trivial involution, then $f$ becomes
a map 
\[
f:X\rightarrow \widehat{\mathcal{F}} 
\]
where $\widehat{\mathcal{F}}\mathcal{\;}$is the space of self-adjoint
Fredholm operators. Now in~\cite{atiyah69:_index_fredh} it is shown that $\widehat{\mathcal{F}}
\mathcal{\;}$has three components 
\[
\widehat{\mathcal{F}}_{+},\widehat{\mathcal{F}}_{-},\widehat{\mathcal{F}}
_{\ast } 
\]
where the first consists of $A$ which are essentially positive (only
finitely many negative eigenvalues), the second is given by essentially
negative operators. These two components are trivial, in the sense that they
are contractible, but the third one is interesting and in fact~\cite{atiyah69:_index_fredh}  
\begin{equation}
\widehat{\mathcal{F}}_{\ast }\sim \Omega \mathcal{F}  \tag{3.2}
\end{equation}
where $\Omega $ denotes the loop space. Since 
\[
\lbrack X,\Omega \mathcal{F}]\cong K^{1}(X) 
\]
this is in agreement with (2.6)\ -- though to get this we have to discard
the two trivial components of $\widehat{\mathcal{F}},\mathcal{\;}$a
technicality to which we now turn.

Lying behind the isomorphism (3.1) is Kuiper's
Theorem~\cite{kuiper65:_hilber} on the contractibility of the unitary
group of Hilbert spaces. Hence to establish that our putative
definition of $K_{\pm }$ coincides with the definition given in
Section 2 we should expect to need a generalization of Kuiper's
Theorem incorporating the involution $A\rightarrow A^{\ast }$ on
operators.  The obvious extension turns out to be false, precisely
because $\widehat{ \mathcal{F}}$, the fixed-point set of *\ on\
$\mathcal{F}$, has the additional contractible components. There are
various ways we can get round this but the simplest and most natural
is to use ``stabilization''. Since $ H\cong H\oplus H$ we can always
stabilize by adding an additional factor of $H.$ In fact Kuiper's
Theorem has two parts in its proof:

\begin{enumerate}
\item[(1)]  The inclusion $U(H)\rightarrow U(H\oplus H)$ defined by $
u\rightarrow u\oplus 1$ is homotopic to the constant map.

\item[(2)]  This inclusion is homotopic to the identity map given by the
isomorphism $H\cong H\oplus H.$
\end{enumerate}

The proof of (1) is an older argument (sometimes called the ``Eilenberg
swindle''), based on a correct use of the fallacious formula 
\begin{eqnarray*}
1 &=&1+(-1+1)+(-1+1)..... \\
&=&(1+-1)+(1+-1)+..... \\
&=&0.
\end{eqnarray*}
The trickier part, and Kuiper's contribution, is the proof of (2).

For many purposes, as in $K$-theory, the stronger version is a luxury and
one can get by with the weaker version (1), which applies rather more
generally. In particular (1) is consistent with taking adjoints (i.e.
inverses in $U(H)$), which is the case we need.

With this background explanation we now introduce formally our second
definition, and to distinguish it temporarily from $K_{\pm }$ as defined in
Section 2, we put 
\begin{equation}
\mathcal{K}_{\pm }(X)=[X,\mathcal{F}]_{\ast }^{s}  \tag{3.3}
\end{equation}
where $\ast $ means we use $Z_{2}$-maps compatible with $\ast $ and $s$ 
means that we use \textbf{stable }homotopy equivalence. More precisely the $
Z_{2}$-maps 
\[
f:X\rightarrow \mathcal{F}(H)\;\;\;\;g:X\rightarrow \mathcal{F}(H) 
\]
are called stably homotopic if the ``stabilized'' maps 
\[
f^{s}:X\rightarrow \mathcal{F}(H\oplus H)\;\;\;\;\;g^{s}:X\rightarrow 
\mathcal{F}(H\oplus H) 
\]
given by $f^{s}=f\oplus J,\;g^{s}=g\oplus J$ are homotopic, where $J$ is a
fixed (essentially unique) automorphism of $H$ with 
\begin{equation}
J=J^{\ast },\;J^{2}=1,\;\;\;+1\;\text{and }-1\;\text{both of infinite
multiplicity}  \tag{3.4}
\end{equation}

Note that under such stabilization the two contractible components $\widehat{
\mathcal{F}}_{+}$ and $\widehat{\mathcal{F}}_{-}$ of $\widehat{\mathcal{F}}
(H)$ both end up in the interesting component $\widehat{\mathcal{F}}_{\ast
}$ of $\widehat{\mathcal{F}}(H\oplus H).$

The first thing we need to observe about $\mathcal{K}_{\pm }(X)$ is that it
is an abelian group. The addition can be defined in the usual way by using
direct sums of Hilbert spaces.   Moreover we can define the negative degree
groups $\mathcal{K}_{\pm }^{-n}(X)$ (for $n\geq 1$) by suspension (with
trivial involution on the extra coordinates), so that 
\[
\mathcal{K}_{\pm }^{-n}(X)=\mathcal{K}_{\pm }(X\times S^{n},X\times \infty
). 
\]
However, at this stage we do not have the periodicity theorem for $\mathcal{K
}_{\pm }(X).$ This will follow in due course after we establish the
equivalence with $K_{\pm }(X).$ As we shall see our construction of (4.2)
is itself closely tied to the periodicity theorem.

Our aim in the subsequent sections will be to show that there is a natural
isomorphism 
\begin{equation}
\mathcal{K}_{\pm }(X)\cong K_{\pm }(X).  \tag{3.5}
\end{equation}
This isomorphism will connect us up naturally with Dirac operators and so
should tie in nicely with the physics.

\section*{4. Construction of the map}

Our first task is to define a natural map 
\begin{equation}
K_{\pm }(X)\rightarrow \mathcal{K}_{\pm }(X).  \tag{4.1}
\end{equation}

We recall from (2.5) that 
\begin{eqnarray*}
K_{\pm }(X) &=&K_{Z_{2}}(X\times R^{1,1}) \\
&=&K_{Z_{2}}(X\times S^{2},\;X\times \infty )
\end{eqnarray*}
where $S^{2}$ is the 2-sphere obtained by compactifying\ $R^{1,1}$, and $
\infty $ is the added point. Note that $Z_{2}$ now acts on $S^{2}$ by a 
\textbf{reflection, }so that it reverses its orientation.

Thus to define a map (4.1) it is sufficient to define a map
\begin{equation}
K_{Z_{2}}(X\times S^{2})\rightarrow \mathcal{K}_{\pm }(X).  \tag{4.2}
\end{equation}
This is where the Dirac operator enters. Recall first that, if we
ignore involutions, there is a basic map
\begin{equation}
K(X\times S^{2})\rightarrow \lbrack X,\mathcal{F}]\cong K(X)  \tag{4.3}
\end{equation}
which is the key to the Bott periodicity theorem.  It is given as
follows.  Let $D$ be the Dirac operator on $S^{2}$ from positive to
negative spinors) and let $V$ be a complex vector bundle on $X\times
S^{2}$, then we can extend, or couple, $D$ to $V$ to get a family
$D_{V}$ of elliptic operators along the $S^{2}$-fibres. Converting
this, in the usual way, to a family of (bounded) Fredholm operators we
get the map (4.3).

We now apply the same construction but keeping track of the involutions. The
new essential feature is that $Z_{2}$ reverses the orientation of $S^{2}$ 
and hence takes the Dirac operator $D$ into its adjoint $D^{\ast }.$ This
is precisely what we need to end up in $\mathcal{K}_{\pm }(X)$ so defining
(4.2).

\begin{rem} Strictly speaking the family $D_{V}$ of
Fredholm operators does not act in a fixed Hilbert space, but in a bundle of
Hilbert spaces. 
The
problem can be dealt with by adding a trivial operator acting on a
complmentary bundle $W$ (so that $W+V$ is trivial).
\end{rem}

% {\em Having defined (4.2) we still need to see that it vanishes on elements
% coming from $K_{Z_{2}}(X)$, i.e. for bundles $V$ that do not depend on the 
% $S^{2}$ factor. But the Dirac operator $D$ on $S^{2}$ is \textbf{
% invertible }and hence (when $V$ comes from $X)$ so are all the operators
% in the family $D_{V}.$ Since we can use the weak form of Kuiper's Theorem
% for the involution *,\ and since we have defined $\mathcal{K}_{\pm }(X)$ by
% stabilization, it follows that an invertible\ family of Fredholm operators
% gives zero in $\mathcal{K}_{\pm }(X).$ Thus (4.2) is well-defined.}

\section*{5. Equivalence of definition}

Let us sum up what we have so far. We have defined a natural homomorphism 
\[
K_{\pm }(X)\rightarrow \mathcal{K}_{\pm }(X) 
\]
and we know that this is an isomorphism for spaces $X$ with trivial
involution --- both groups coinciding with $K^{1}(X)$.  Moreover, if
for general $X$, we ignore the involutions, or equivalently replace
$X$ by $X\times \{0,1\}$, we also get an isomorphism, both groups now
coinciding with $K^{0}(X)$.

General theory then implies that we have an isomorphism for all $X.$ We
shall review this general argument.

Let $A,B$ be representable theories, defined on the category of $Z_{2}$
-spaces, so that 
\begin{eqnarray*}
A(X) &=&[X,\mathcal{A}] \\
B(X) &=&[X,\mathcal{B}]
\end{eqnarray*}
where $[\phantom{x},\phantom{x}]$ denotes homotopy classes of
$Z_{2}$-maps into the classifying spaces of $\mathcal{A},\mathcal{B}$
of the theories. A natural map $ A(X)\rightarrow B(X)$ then
corresponds to a $Z_{2}$-map\ $\mathcal{ A\rightarrow B}$. Showing
that $A$ and $B$ are isomorphic theories is equivalent to showing that
$\mathcal{A}$ and $\mathcal{B}$ are $Z_{2}$-homotopy equivalent.

If we forget about the involutions then isomorphism of theories is the same
as ordinary homotopy equivalence. Restricting to spaces $X$ with trivial
involution corresponds to restricting to the fixed-point sets of the
involution on $\mathcal{A\;}$and $\mathcal{B}.$

Now there is a general theorem in homotopy theory~\cite{james78} which
asserts (for reasonable spaces including Banach manifolds such as
$\mathcal{F}$) that, if a $Z_{2}$-map $\mathcal{A\rightarrow B}$\ is
both a homotopy equivalence ignoring\ the involution and for the
fixed-point sets, then it is a $Z_{2}$-homotopy
equivalence. Translated back into the theories $A,B$ it says that the
map $A(X)\rightarrow B(X)$ is an isomorphism provided it holds for
spaces $X$ with the trivial $Z_{2}$-action, and for $Z_{2}$-spaces $X$ of
the form $Y\times \{0,1\}$.

This is essentially the situation we have here with 
\[
A=K_{\pm}\qquad B=\mathcal{K}_{\pm}.
\]
Both are representable. The representability of the first 
\[
K_{\pm }(X)\cong K_{Z_{2}}(X\times R^{1,1}) 
\]
arises from the general representability of $K_{Z_{2}}$, the classifying\
space being essentially the double loop space of $\mathcal{F}(H\otimes
C^{2})$ with an appropriate $Z_{2}$-action. The second is representable
because 
\begin{equation}
\mathcal{K}_{\pm }(X)=[X,\mathcal{F}]_{\ast }^{s}=[X,\mathcal{F}_{s}]_{\ast }
\tag{5.1}
\end{equation}
where $\mathcal{F}^{s}$ is obtained by stabilizing $\mathcal{F}.$ More
precisely 
\[
\mathcal{F}^{s}=\lim_{n\rightarrow \infty }\mathcal{F}_{n} 
\]
where $\mathcal{F}_{n}=\mathcal{F}(H\otimes C^{n})$ and the limit is taken
with respect to the natural inclusions, using $J$ of (3.4)\ as a base
point. The assertion in (5.1)\ is easily checked and it simply gives two
ways of looking at the stabilization process.

We have thus established the equivalence of our two definitions $K_{\pm }$ 
and $\mathcal{K}_{\pm }.$

\section*{6. Free involutions}

We shall now look in more detail at the case of free involutions and,
following the notation of Section 1, we shall denote the free $Z_{2}$-space
by $\tilde{X}$ and its quotient by $X.$

The reason for introducing the stabilization process in Section 3 concerned
fixed points. We shall now show that, for free involutions, we can dispense
with stabilization. Let 
\[
\mathcal{F}\rightarrow \mathcal{F}^{s} 
\]
be the natural inclusion of $\mathcal{F}$ in the direct limit space. This
inclusion is a $Z_{2}$-map and a homotopy equivalence, though \textbf{not }a 
$Z_{2}$-homotopy equivalence (because of the fixed points). Now given the
double covering $\tilde{X}\rightarrow X$ we can form the associated fibre
bundles $\mathcal{F}_{X}$ and $\mathcal{F}_{X}^{s}$ over $X$ with fibres $
\mathcal{F\;}$and $\mathcal{F}^{s}.$ Thus 
\[
\mathcal{F}_{X}=\tilde{X}\times_{Z_{2}}\mathcal{F\;\;\;\;F}_{X}^{s}=\tilde{X
}\times_{Z_{2}}\mathcal{F}^{s} 
\]
and we have an inclusion 
\[
\mathcal{F}_{X}\rightarrow \mathcal{F}_{X}^{s} 
\]
which is fibre-preserving. This map is a homotopy-equivalence on the
fibres and hence, by a general theorem~\cite{dold63:_partit} (valid in
particular for Banach manifolds) a fibre homotopy equivalence.   It
follows that the homotopy classes of sections of these two fibrations
are isomorphic. But these are the same as
\[
\left[ \tilde{X},\mathcal{F}\right]_{\ast }\;\text{and\ }\left[ \tilde{X},
\mathcal{F}\right]_{\ast }^{s}=\mathcal{K}_{\pm }(\tilde{X}). 
\]
This show that, for a free involution, we can use $\mathcal{F\;}$instead of $
\mathcal{F}^{\text{s}}.$ Moreover it gives the following simple description
of $K_{\pm }(\tilde{X}):$
\begin{equation}
K_{\pm }(\tilde{X})=\;\text{Homotopy classes of sections of }\mathcal{F}
_{X}\;.  \tag{6.1}
\end{equation}
This is the $K$-theory analogue of twisted cohomology described in Section
1. A corresponding approach to the higher twist of $K$-theory given by an
element of $H^{3}(X;Z)$ will be developed in~\cite{atiyah:_twist_k}

\section*{7. The real case}

Everything we have done so far extends, with appropriate modifications, to
real $K$-theory. The important difference is that the periodicity is now 8
rather than 2 and that, correspondingly we have to distinguish carefully
between \textbf{self-adjoint }and \textbf{skew-adjoint\ }Fredholm operators.
Over the complex numbers multiplication by $i$ converts one into the other,
but over the real numbers there are substantial differences.

We denote by $\mathcal{F}^{1}(R)$ the interesting component of the space of
real self-adjoint Fredholm operators $\widehat{\mathcal{F}}(R)$ on a real
Hilbert space (discarding two contractible components as before). We also
denote by $\mathcal{F}^{-1}(R)$ the space of skew-adjoint Fredholm
operators. Then in~\cite{atiyah69:_index_fredh} it is proved that 
\begin{align*}
[X,\mathcal{F}^{1}(R)] &\cong KO^{1}(X)  \tag{7.1} \\
\lbrack X,\mathcal{F}^{-1}(R)] &\cong KO^{-1}(X)\cong KO^{7}(X) 
\tag{7.2}
\end{align*}
showing that these are essentially different groups.

Using (7.1), stabilizing, and arguing precisely as before, we define 
\begin{align*}
KO_{\pm }(X) &= KO_{Z_{2}}^{1}(X\times R^{1,0}) 
\cong KO_{Z_{2}}(X\times R^{1,7})\\
\mathcal{KO}_{\pm }(X) &=[X,\mathcal{F}(R)]_{\ast }^{s},
\end{align*}
where (in (2.5)) the mod $2$ periodicity of K has been replaced by the
mod $8$ periodicity of $KO$.
But we cannot now just use the Dirac operator on $S^{2}$ because this is
not real.  Instead we have to use the Dirac operator on $S^{8},$ 
% This means
% we should suspend 6 times and use 
% \[
% KO_{\pm }^{2}(X)=KO_{\pm }^{-6}(X)=KO_{Z_{2}}(X\times R^{1,7})\;.
% \]
% The Dirac operator on $S^{8}$ 
which then gives us our map 
\begin{equation}
KO_{\pm }(X)\rightarrow \lbrack X,\mathcal{F}(R)]_{\ast }^{s}\;. 
\tag{7.3}
\end{equation}
% \begin{equation}
% KO_{\pm }^{2}(X)\rightarrow \lbrack X,\mathcal{F}(R)]_{\ast }^{s}\;. 
% \tag{7.3}
% \end{equation}
The same proof as before establishes the isomorphism of (7.3), so that 
% \[
% KO_{\pm }^{2}(X)\cong \mathcal{KO}_{\pm }(X)
% \]
\[
KO_{\pm }(X)\cong \mathcal{KO}_{\pm }(X)
\]
and more generally for $q$ modulo 8, 
\begin{equation}
KO_{\pm }^{q}(X)\cong \mathcal{KO}_{\pm }^{q}(X).  \tag{7.4}
\end{equation}
% \begin{equation}
% KO_{\pm }^{q}(X)\cong \mathcal{KO}_{\pm }^{q-2}(X).  \tag{7.4}
% \end{equation}

In~\cite{atiyah69:_index_fredh} there is a more systematic analysis of
Fredholm operators in relation to Clifford algebras and using this it
is possible to give more explicit descriptions of $KO_{\pm }^{q}(X)$,
for all $q$, in terms of $Z_{2}$ -mappings into appropriate spaces of
Fredholm operators. This would fit in with the different behaviour of
the Dirac operator in different dimensions (modulo 8).

\bibliographystyle{amsplain}
\providecommand{\bysame}{\leavevmode\hbox to3em{\hrulefill}\thinspace}
\providecommand{\MR}{\relax\ifhmode\unskip\space\fi MR }
% \MRhref is called by the amsart/book/proc definition of \MR.
\providecommand{\MRhref}[2]{%
  \href{http://www.ams.org/mathscinet-getitem?mr=#1}{#2}
}
\providecommand{\href}[2]{#2}

\end{document}